\documentclass[11pt]{article}
\usepackage{times}
\usepackage{indentfirst}
\usepackage{amsthm,amsmath, amssymb}
\usepackage{calrsfs}
\usepackage{graphicx}
\usepackage{multicol}
\usepackage{caption}
\usepackage{subcaption}
\usepackage{xcolor}

\usepackage[mathscr]{euscript}

\usepackage{algorithm}
\usepackage[noend]{algpseudocode}

\makeatletter
\def\BState{\State\hskip-\ALG@thistlm}
\makeatother

\newtheorem{Theorem}{Theorem}
\newtheorem{Lemma}{Lemma}
\newtheorem{Proposition}{Proposition}

\newtheorem{Remark}{Remark}
\newtheorem{Example}{Example}

\newcommand{\eq}[1]{(\ref{eq:#1})}

\newcommand{\thr}[1]{Theorem~\ref{thr:#1}}

\newenvironment{Proof} 
	{\par\noindent{\bf Proof.}} 
	{\hfill$\scriptstyle\blacksquare$} 

\newcommand{\PB}{\mathbb{P}}
\newcommand{\Expect}{\mathbb{E}}

\begin{document}

\title{\bf \Large Moments of the first descending epoch \\
for a random walk with negative drift}

\author{Sergey Foss,\thanks{Postal adress: School of Mathematical and Computer Sciences, Heriot-Watt University,
EH14 4AS, Edinburgh, United Kingdom. Email address:
sergueiorfoss25@gmail.com.} {\em
\normalsize Heriot-Watt University} 
\and Timofei Prasolov,
\thanks{Postal address:  Novosibirsk State University, 630090, Novosibirsk, Russia. Email
address:prasolov.tv@yandex.ru} {\em \normalsize Novosibirsk State University}}
\date{\today}



\maketitle

\begin{abstract}
We consider the first descending ladder epoch $\tau = \min \{n\ge 1 : S_n\le 0\}$ of a random walk $S_n = \sum_1^n \xi_i, n\ge 1$ with i.d.d. summands having a negative drift 
${\mathbb E} \xi = -a< 0$. Let $\xi^+ = \max (0, \xi_1)$.
It is well-known that, for any $\alpha>1$, the finiteness of ${\mathbb E}(\xi^+)^{\alpha}$ implies the finiteness of ${\mathbb E} \tau^{\alpha}$ and, for any $\lambda >0$, the finiteness of ${\mathbb E} \exp({\lambda\xi^+})$ implies that of ${\mathbb E} \exp({c\tau})$ where $c>0$ is, in general, another constant that depends 
on the distribution of $\xi_1$. We consider the intermediate case, assuming that 
${\mathbb E} \exp({g(\xi^+)})<\infty$ for a positive increasing function $g$ such that $\liminf_{x\to\infty} g(x)/\log x = \infty$ and $\limsup_{x\to\infty}  g(x)/x =0$, and that ${\mathbb E} \exp({\lambda\xi^+})=\infty$, for all $\lambda>0$. Assuming a few further technical assumptions, we show that then    ${\mathbb E} \exp({(1-\varepsilon){g}((1-\delta)a\tau)})<\infty$, for any $\varepsilon,\delta \in (0,1)$.

\end{abstract}

\begin{quotation}
\end{quotation}

\section{Introduction and the main result}

Let $\xi, \xi_1, \xi_2, \ldots, \xi_n, \ldots$ be independent and identically distributed (i.i.d.) random variables (r.v.'s) with a common distribution function $F$ having a finite negative  mean $\Expect \xi  = -a < 0$. Let $S_0 = 0$, $S_n = \sum_{k=1}^n \xi_k, n\ge 1$ be a  random walk, and $\tau = \min\{n\ge 1: \ S_n \le 0\}<\infty$ a.s. its first descending ladder epoch.

The descending ladder epoch  $\tau$ plays an important role
in theoretical and applied probability. 
In particular, 
 $\tau$  represents the length of a busy cycle in a $GI/GI/1$ queueing system. Namely, consider a FIFO single-server queue with i.i.d. interarrival times $\{t_n\}$ with a finite mean ${\mathbb E} t_1=a$ and independent of them i.i.d. service times $\{\sigma_n\}$ with a finite mean ${\mathbb E} \sigma_1=b<a$.
 Let $W_n$ be the waiting time of customer $n$. Assume $W_1=0$, i.e. customer 1 arrives at an empty queue. The sequence $\{W_n\}$ satisfies the {\it Lindley recursion}
 \begin{align}\label{waiting}
 W_{n+1} = \max (0, W_n+\sigma_n-t_n) \quad n\ge 1.
 \end{align}
 We may let $\xi_n=\sigma_n-t_n$ and conclude that
 $\tau$ is the number of customers served in the queue during the {\it first busy cycle}, i.e. customer $\tau+1$ is the next customer after customer 1 that finds the queue empty.
 
 We are interested in the existence (finiteness) of moments of $\tau$ in terms of moments of the common distribution $F$ of the summands.   In particular, the existence of a power (or an exponential) moment of $\tau$ implies corresponding convergence rates in stability and continuity theorems for various single- and multi-server queueing systems, see e.g. Theorems 2 and 11 in Chapter 4 of \cite{Bor}.

The following results are known
(see, e.g., Theorems III.3.1 and 3.2 in  \cite{Gut}, and also \cite{Hey}). Let $\alpha >1$ and $\lambda>0$.
\begin{align}\label{power}
 \mbox{If} \ \ \Expect (\xi^+)^\alpha<\infty, \ \ \mbox{then}
\ \  \Expect \tau^\alpha < \infty.
\end{align}
\begin{align}\label{exp}
 \mbox{If} \ \ 
\Expect \exp(\lambda\xi)<\infty,   
\ \ \mbox{then there exists} \ \ 
c >0 \ \ \mbox{(that depends on} \ \ F) \ \ \mbox{such that} \ \  \Expect \exp(c\tau)<\infty.
\end{align}

One can view \eqref{power} and \eqref{exp} as two particular cases of the following
implication:
\begin{align}\label{general}
\mbox{If} \ \ {\mathbb E} G(\xi^+)<\infty, \ \ \mbox{then} \ \  
{\mathbb E} G(C\tau)<\infty, \ \ \mbox{for a certain}
\ \  C>0.
\end{align} 
Indeed, \eqref{power} is a particular case of \eqref{general} with $G(x)=x^{\alpha}$, and \eqref{exp} a particular case of \eqref{general}
with $G(x)=\exp (\lambda x)$ (clearly, for $\lambda >0$, exponential moments  $\Expect \exp(\lambda\xi)$ and 
$\Expect \exp(\lambda\xi^+)$ are either finite or infinite simultaneously).

In this article, we consider the intermediate case where $G$ is a monotone function that increases faster than any power function and slower than any exponential function. It is convenient to us to use representation $G(x) = e^{g(x)}$ and work with function $g$ instead.
Here is our main result. 

\begin{Theorem}\label{thr:main}
Assume that ${\mathbb E} \exp (c\xi)=\infty$, for any $c>0$. If  a function $g$ satisfies conditions $(C1)-(C3)$, introduced below, and if 
\begin{align}\label{eq:ex1}
\Expect \exp (g(\xi)) < \infty,
\end{align} 
then 
\begin{align}\label{eq:ex2}
\Expect \exp ((1-\varepsilon)g((a-\delta)\tau)) < \infty, 
\ \ 
 \mbox{for any} \ \  \varepsilon \in (0, 1) \ \ \mbox{and} \ \  \delta \in (0, a).
 \end{align}
\end{Theorem}

The conditions $(C1)-(C3)$ are as follows: 
\begin{itemize}
\item{(C1)} function $g$ is positive, increasing and differentiable;
\item{(C2)} $\lim_{x \to \infty} g'(x) = 0$;
\item{(C3)} there exist a constant $\gamma \in (0, 1)$ such that
\begin{align}\label{eq:C3_3}
\int_{1}^\infty \exp(-(1-\gamma)g(x))dx < \infty
\end{align}
 and positive constants $x_0$ and $A$ such that, for any $x_0 < y \le x/2$, 
\begin{align}\label{eq:C3_2}
g(x) - g(x-y) \le \gamma g(y) + A.
\end{align}
\end{itemize}

It follows from condition $(C2)$ that $\sup_{x\ge x_0} g'(x) \downarrow 0$ as $x_0\to\infty$. Therefore, we may choose $x_0$ in condition $(C3)$ and constant $B>0$ such that
\begin{align}\label{eq:C3_1}
  g'(x) < B, \ \ \text{for $x > x_0$}.
\end{align}

\begin{Remark}
\normalfont
Conditions $(C1)-(C3)$ are given in the form that are convenient to us, they may be weakened. 
For example, it is not necessary to assume differentiability,
and condition $(C2)$ can be adjusted to `dying'  growth rate that also gives us  inequality \eq{C3_1}. However, inequalities \eq{C3_3} and \eq{C3_2}  are more substantial since they target heavy-tailed ``Weibull-type'' and ``lognormal-type'' distributions.
\end{Remark}

\begin{Example}
\normalfont
Here are examples of functions $g$ that satisfy conditions $(C1)-(C3)$: \\
$g_1(x)= (\log \max (x,1))^{\alpha}$, 
 $g_2(x) = (x^+)^\beta$ and $g_3(x)= (x^+)^{\beta} \log (max (x,1))$, where $\alpha >1$ and $\beta \in (0, 1)$. More generally, the functions $g_1$ and $g_3$  continue to satisfy condition $(C1)-(C3)$ if the logarithmic function therein is replaced by  a ``sufficiently smooth'' increasing and slowly varying function.
\end{Example}

\begin{Remark}
\normalfont
Note that one can represent  \eqref{eq:ex2} in an equivalent form as:
  \begin{align*}
\Expect \exp ((1-\varepsilon)g((1-\varepsilon)a\tau)) < \infty, 
\ \ 
 \mbox{for any} \ \  \varepsilon \in (0, 1).
 \end{align*}
On the other hand, given condition  \eq{ex2}, the inequality in \eq{ex2} also holds for function $g_1$ from Example 1 with $\delta = 0$ and any $\varepsilon \in (0,1)$, and for functions $g_2$ and $g_3$ with $\varepsilon =0$
and any $\delta\in (0,a)$. 
Let us show this for $g_3$. Indeed, for any $\delta_1\in (0,a)$ there exist  $\varepsilon_2\in (0,1)$ and $\delta_2\in (0,a)$ such that 
$$\lim_{x \to \infty}\frac{g_3((a-\delta_1)x)}{(1-\varepsilon_2)g_3((a-\delta_2)x)} = \frac{(a-\delta_1)^\beta}{(1-\varepsilon_2)(a-\delta_2)^\beta} < 1.$$
Then there exists a constant $c>0$ such that
$$\Expect \exp (g_3((a-\delta_1)\tau)) \le c\Expect \exp ((1-\varepsilon_2)g_3((a-\delta_2)\tau)) < \infty.$$
\end{Remark} 

Our proof of the theorem includes two steps. First, we show the existence of a r.v. $\widetilde{\xi\,} \ge_{st}\xi$ that has a strong subexponential distribution, negative mean and certain finite moments. Second, we prove that the stopping time for the random walk with new increments $\{{\widetilde{\xi}}_n\}$ satisfies the conditions of the theorem.

We use the following notation and conventions. For a distribution function $F$ on the real line, $\overline{F}(x)=1-F(x)$ is its {\it tail distribution} function. For two strictly positive functions $h_1$ and $h_2$, equivalence $h_1(x)\sim h_2(x)$ means that $\lim_{x\to\infty} h_1(x)/h_2(x)=1$. For two r.v.s $\eta_1$ and $\eta_2$,
stochastic inequality $\eta_1\le_{st}\eta_2$ means that ${\mathbb P} (\eta_1>x) \le {\mathbb P} (\eta_2>x)$, for all $x$. For an increasing function $g$, its (generalised) {\it inverse function} $g^{-1}$ is defined as $g^{-1}(t) = \inf\{x:  g(x) > t\}$. Then the sets $\{g(x)>t\}$ and $\{x>g^{-1}(t)\}$ do coincide. A function $f$ is {\it slowly varying} if $f(\lambda x)/x \to 1$, as $x\to\infty$, for $\lambda > 0$, and {\it regularly varying with exponent} $\alpha$  if $f(\lambda x)/x \to \lambda^\alpha$.

\section{Proof of the theorem}\label{sect:proofs}

Recall the following definitions. 
Let $ F $ be a distribution on the real line with right-unbounded support.
We say that $F$ is \emph{long-tailed} if $\lim_{x \to \infty} \overline{F}(x-1) / \overline{F}(x) = 1$. Since the tail function $\overline{F}$ is monotone non-increasing, its long-tailedness implies that $\lim_{x \to \infty} \overline{F}(x-y) / \overline{F}(x) = 1$, for any $y>0$.

Further, let a distribution $F$ have  right-unbounded support and finite mean $m=\int_0^{\infty} \overline{F}(y) dy$ on the positive half line. We say that $F$ is \emph{strong subexponential} and write $F \in \mathscr{S}^*$ if $\int\limits_0^x \overline{F}(x-y)\overline{F}(y) dy \sim 2m \overline{F}(x)$, as $x \to \infty$. 
The strong subexponentiality is a {\it tail property}:
if a distribution function $F$ is  strong subexponential and if $G$ is another distribution function such that $\overline{F}(x)\sim \overline{G}(x)$, then $G$ is also strong subexponential (see, e.g., \cite{FosKorZac}, Theorem 3.11).

\subsection{Step one: an upper-bound random variable having a strong subexponential distribution}

Let $\zeta = \exp(g(\xi))$. Since $\Expect \zeta < \infty$, $\PB\{\zeta > x\} = o(x^{-1})$, as $x\to\infty$. Then, in particular, one can choose $K \ge \exp(g(x_0))$ (where constant $x_0$ is from condition $(C3)$) such that $\PB\{\zeta > x\}\le Kx^{-1}$, for all $x>0$, and 
introduce a new non-negative  r.v. $\widehat{\zeta}$  with the tail distribution 
\begin{equation}\label{eq:Kx}
\PB\{\widehat{\zeta} > x\} = \min(1, Kx^{-1}), \ x\ge 0.
\end{equation}
Clearly, $\PB\{\zeta > x\} \le \PB\{\widehat{\zeta} > x\}$, for all $x$, and $\Expect \widehat{\zeta}^{1-\varepsilon} < \infty$, for $\varepsilon \in (0, 1)$.

\begin{Lemma} Under the assumptions $(C1)-(C3)$, 
the r.v. $\widehat{\xi} = g^{-1}(\ln(\widehat{\zeta}))$ has a strong subexponential distribution $\widehat{F}$.
\end{Lemma}

\begin{Proof}
We use the following result.
\begin{Proposition}\label{thr:criteria for subexponality}
(This is a part of Theorem 3.30 from \cite{FosKorZac}). Let $F$ be a long-tailed distribution on the real line.
Let $R(x) = -\ln \overline{F}(x)$. Suppose that there exist $\gamma < 1$ and $A'<\infty$ such that
\begin{equation}\label{eq: condition of less than linear function g}
R(x) - R(x-y) \le \gamma R(y) + A',
\end{equation}
for all $x>0$ and $y\in [0, x/2]$. If, in addition, 
\begin{align}\label{eq:int}
\mbox{the function} \ \ \exp(-(1-\gamma)R(x)) \ \ 
\mbox{ is integrable over} \ \  [0,\infty),
\end{align}
 then $F\in\mathscr{S}^*$.
\end{Proposition}

To apply Proposition \ref{thr:criteria for subexponality},  we need to verify the
long-tailedness of $\widehat{F}$ and conditions \eqref{eq: condition of less than linear function g} and \eqref{eq:int}.
First, we show the long-tailedness of $\widehat{F}$. For a fixed $y>0$ and large $x$, we have $\PB\{\widehat{\xi} > x+y\} = K\exp(-g(x+y))$.

From the first-order Taylor expansion $g(x + y) = g(x) + y g'(z)$, for some $z\in (x, x+y)$, and from condition $(C2)$ we get
\begin{equation*}
1 
\ge \frac{\PB\{\widehat{\xi} > x+y\}}{\PB\{\widehat{\xi} > x\}} \ge \frac{\exp(-g(x) - y g'(z))}{\exp(-g(x))} = \exp(-y g'(z)) = \exp(o(1)) = (1 + o(1)),
\end{equation*}
as $x\to \infty$. Thus, the distribution of $\widehat{\xi}$ is long-tailed.

Second, we verify condition \eqref{eq: condition of less than linear function g}.  It is equivalent to 
\begin{equation}\label{eq:transform of less-linear condition}
\frac{\overline{H}(x-y)}{\overline{H}(x)} \le \frac{\exp(A')}{\overline{H}^\gamma(y)},
\end{equation}
where $\overline{H}(x) = \PB\{ \widehat{\zeta} >  \exp(g(x))\}$. We take $\gamma$ from condition $(C3)$. Next we show the existence of an appropriate constant $A'$.  

Let $x_1 = \inf\{x: \ \overline{H}(x) < 1\}$. Since we have chosen $K > \exp(g(x_0))$, we get $x_1 \ge x_0$. We consider four cases depending on whether $\overline{H}(x) = 1$ or $\overline{H}(x) = K\exp(-g(x))$. 

Assume $x \le x_1$. Then inequality \eqref{eq:transform of less-linear condition} holds if  we take $A'\ge0$.

Assume $x -y \le x_1 < x$. Then \eqref{eq:transform of less-linear condition} is equivalent to $K^{-1}\exp(g(x)) \le \exp(A')$. Since $x/2 \le x-y \le x_1$, inequality $A'\ge g(2x_1) - \ln K$ is a sufficient condition on $A'$ to satisfy \eqref{eq:transform of less-linear condition}.

Assume $y \le x_1 < x-y$. Then \eqref{eq:transform of less-linear condition} is equivalent to $\exp(g(x) - g(x-y)) \le \exp(A')$. Since $g(x-y)  = g(x) - y g'(z)$, for $z \in (x-y, x)$, we have $g(x) - g(x-y) = y g'(z) < B x_1$. Therefore, it is sufficient to assume $A' \ge B x_1$.

Next, assume $y > x_1$. Then \eqref{eq:transform of less-linear condition} is equivalent to $\exp(g(x) - g(x-y)) \le K^{-\gamma}\exp(\gamma g(y) + A')$. From condition $(C3)$ it is sufficient to assume $A' \ge A + \gamma \ln K$ for the Proposition 2 to hold.

Finally, condition \eqref{eq:int} follows directly from $(C3)$.
\end{Proof}

By construction, $ \Expect \exp((1-\varepsilon)g(\widehat{\xi})) = \Expect \widehat{\zeta}^{1-\varepsilon} < \infty$. However, we need our upper-bound to have sufficiently close mean to the original. Thus, we need the following lemma.

\begin{Lemma}\label{thr:transformation to have negative mean}
Assume that conditions $(C1)-(C3)$ hold. For any $\delta \in (0, a)$, we can introduce a r.v. $\widetilde{\zeta}$ such that $\widetilde{\xi} = g^{-1}(\ln(\widetilde{\zeta}))$
has a strong subexponential distribution, $\widetilde{\xi}\ge_{st} \xi$  and, in addition, $\Expect \widetilde{\xi} < \Expect \xi + \delta = -a + \delta < 0.$
\end{Lemma}
\begin{Proof}
Since the distributions of $\xi$ and $\widehat{\xi}$ have right-unbounded support, for all $V>0$ we can find $V' > V$ such that there exists r.v. $\widetilde{\xi} $ with right tail
\begin{equation*}
\PB\{\widetilde{\xi} > t\} = \begin{cases}
\PB\{\xi > t\}, & t < V,\\
\PB\{\xi > V\}, & V \le t < V',\\
\PB\{\widehat{\xi} > t\}, & t\ge V'.
\end{cases}
\end{equation*}
Clearly, $\xi \le_{st} \widetilde{\xi} \le_{st} \widehat{\xi}$. Since $\widetilde{\xi}$ and $\widehat{\xi}$ have the same right tail, $\widetilde{\xi}$ has a strong subexponential distribution. By choosing sufficiently large $V$ we can make $\Expect \widetilde{\xi} =  \Expect\left(  \widetilde{\xi}; \ \widetilde{\xi} \le V \right) + \int_V^\infty\PB\{\widetilde{\xi} > t\} dt  = \Expect\left(  \xi; \ \xi \le V \right) + \int_V^\infty\PB\{\widetilde{\xi}> t\}dt < -a+\delta$.
\end{Proof}

\subsection{Step two: existence of moments of the first descending epoch for strong subexponential distributions}

We have introduced a r.v. $ \widetilde{\xi}$ with negative drift $\Expect  \widetilde{\xi} =  -\widetilde{a} < -a+\delta < 0$ and a finite moment\\ $\Expect \exp((1-\varepsilon)g( \widetilde{\xi})) < \infty$, such that $\xi \le_{st} \widetilde{\xi}$. Now we want to show that the stopping time $ \widetilde{\tau}$ satisfies\\ $\Expect \exp((1-\varepsilon)g( (a-\delta) \widetilde{\tau})) < \infty$. 

Without loss of generality, we may assume that the  distribution $\widetilde{F}$ of the r.v.'s $\widetilde{\xi}_k$ is bounded below, i.e. $\widetilde{\xi}_k \ge -L$ a.s., for some $L \in (0, \infty)$.
 Indeed, let us  choose an arbitrary $L>0$ and take $\xi'_i = \max(\widetilde{\xi}_i, -L)$, $i\ge 1$. Then the random walk $S'_0 = 0$, $S'_n = \sum_{k=1}^n \xi'_k$ satisfies $S'_n\ge \widetilde{S}_n$ a.s., for all $n$ and, therefore, $\tau' = \inf\{n\ge 1: \ S'_n \le 0\} \ge \widetilde{\tau}$ a.s.

By taking $L$ large enough, we can make $\Expect \xi' = 
{\mathbb E} \widetilde{\xi} - {\mathbf E} (\widetilde{\xi}+L; \widetilde{\xi} \le -L) $ as close to ${\mathbb E}\widetilde{\xi}$ as one wishes and, in particular, smaller than zero. Since $\sup_{x\le 0} g(x)<\infty$, condition \eqref{eq:ex1} implies the finiteness of $\Expect \exp (g(\xi'))$ too. If we prove the statement of \thr{main} for the random walk with increments $\xi'_n$, then we prove it for the initial random walk, too.

We write $h(\cdot) = (1-\varepsilon)g(\cdot)$ for short. We prove now that $\Expect \exp(h((a - \delta)\widetilde{\tau})) < \infty$. Let $\chi = S_{\widetilde{\tau}}$,  $\chi \in [-L, 0]$. We have
\begin{equation*}
(a - \delta)\widetilde{\tau}= (a - \delta)\widetilde{\tau} + \chi-\chi \leq ((a - \delta)\widetilde{\tau}+\chi)+ L = \sum_{i=1}^{\widetilde{\tau}} (\widetilde{\xi}_{i}+a - \delta)+L. 
\end{equation*}
Let $\psi_i = \widetilde{\xi}_{i}+a - \delta$. Thus, $\Expect \psi_1 < 0$ and, since $\PB\{\psi_1 > x\} \sim \PB\{\widetilde{\xi}_1 > x\}$, r.v. $\psi_1$ has a strong subexponential distribution. 
 From inequality \eq{C3_1} and the first-order Taylor expansion for $h$ we get $h(x + y) \le h(x) + (1-\varepsilon)B y$, for $x> x_0$, and thus, 
\begin{equation*}
\Expect \exp(h((a - \delta)\widetilde{\tau})) \le \Expect \exp\left( h\left(  \sum_{i=1}^{\widetilde{\tau}}\psi_{i} + L\right) \right)
 \le \exp\left( h\left(x_0 + L\right) \right) +  \exp\left( (1-\varepsilon)B L\right) \Expect \exp\left( h\left(  \sum_{i=1}^{\widetilde{\tau}}\psi_{i}\right) \right).
\end{equation*}
Further,
\begin{equation}\label{eq:eq15}
\Expect \exp\left( h\left(  \sum_{i=1}^{\widetilde{\tau}}\psi_{i}\right) \right) = \int_0^\infty \PB\left\lbrace \exp\left( h\left(  \sum_{i=1}^{\widetilde{\tau}}\psi_{i}\right) \right)  > t\right\rbrace  dt \le \int_0^\infty \PB\left\lbrace \sum_{i=1}^{\widetilde{\tau}}\psi_{i} > h^{-1}(\ln t)\right\rbrace dt.
\end{equation}

Next, we need the following result:
\begin{Proposition}\label{th: Upper bound for convolution tail}
(Theorem 1 in \cite{FosZac}). Let $\Expect \psi < 0$ and let $\tau$ be a stopping time for $\{\psi_n\}$. Denote $M_\tau = \max_{0\leq j \leq \tau}\sum_{i=1}^{j}\psi_{i}$ and let $F_{\psi}(x)$ be the distribution function of $\psi$. Under condition $F_\psi \in \mathscr{S}^*$ we have
\begin{equation*}
\lim_{x\to \infty} \frac{P\{M_\tau > x\}}{\overline{F}_\psi(x)} = \Expect\tau.
\end{equation*}
\end{Proposition}

Clearly, we can apply Proposition \ref{th: Upper bound for convolution tail}:
we have $\Expect \psi < 0$, and $\psi$ has a strong subexponential distribution. Also, the r.v. $\widetilde{\tau}$ is a stopping time w.r.t. $\{\widetilde{\xi}_n\}$ and, therefore, w.r.t. $\{\psi_n\}$. Thus, the conditions of Proposition \ref{th: Upper bound for convolution tail} hold. Now, combining this with \eq{eq15}, we get that, for every $\Delta >0$, there exists a constant $N$ such that
\begin{equation*}
\Expect \exp\left( h\left(  \sum_{i=1}^{\widetilde{\tau}}\psi_{i}\right) \right) \le N + \int_N^\infty \PB\{M_{\widetilde{\tau}} > h^{-1}(\ln t)\}dt \le  N + (\Expect\widetilde{\tau} +\Delta )\int_N^\infty \PB\{\widetilde{\xi} + a-\delta > h^{-1}(\ln t)\}dt,
\end{equation*}
and the integral on the right-hand side of the latter inequality is finite. This  concludes the proof of the theorem.

\section{Further comments}

In our theorem, the coefficients $(1-\varepsilon)$ and $(1-\delta)$  
appear because the first moment of the upper-bound distribution in \eqref{eq:Kx} is infinite. 
The following nice result may help to eliminate the coefficients under certain assumptions discussed below.

\begin{Proposition}\label{prop4}
({\em  Corollary 1 in \cite{Den}}) Let $\zeta$ be a nonnegative r.v. and $\Expect \zeta^\alpha < \infty$ for some $\alpha >0$. Then there exists a r.v. $\widehat{\zeta}$ such that $\Expect \widehat{\zeta}^\alpha < \infty$, $\PB\{\widehat{\zeta} > t\}$ is a function of regular variation with exponent $-\alpha$, and  $\zeta \le_{st} \widehat{\zeta}$.
\end{Proposition}

We can apply Proposition \ref{prop4} with $\alpha = 1$, $\zeta = \exp(g(\xi))$, and then the upper bound $\widehat{\zeta}$ has the tail distribution $\PB\{\widehat{\zeta} \ge x\} \sim l(x)/x$, which is integrable. Here $l(x)$ is a slowly varying function.  If in addition $l(x)$ is sufficiently smooth (to be justified), there is a chance to show that $\widehat{\xi} = g^{-1}(\ln \widehat{\zeta})$ has a strong subexponential distribution and $\Expect \exp(g(\widehat{\xi})) < \infty$. Then the statement of the theorem holds with $\varepsilon = \delta =0$.

Another way to apply Proposition \ref{prop4} is to provide an alternative proof of Theorem III.3.1 in  \cite{Gut}. Indeed, in this case the distribution of $(\xi^+)^{\alpha}$ possesses an integrable majorant having a regularly varying distribution. Since any power of a regularly varying function is also a regularly varying function, the distribution of $\xi^+$ possesses a majorant having a regularly varying distribution with finite moment of
order $\alpha$. And it is known that any regularly varying distribution with finite mean is strong subexponential.

{\bf Acknowledgment}. he work is supported by Mathematical Center in Akademgorodok under agreement No. 075-15-2022-282 with the Ministry of Science and Higher Education of the Russian Federation.

\end{document}